\documentclass[12pt]{amsart}
\usepackage{amsmath}
\usepackage{amssymb}
\usepackage{amscd} 
\usepackage{graphicx}

\newcommand{\bz}{\mathbf Z}

\newcommand{\bc}{\mathbf C}

\newcommand{\ind}{\operatorname{ind}}
 
\newcommand{\dir}{\mathcal  D}
\newcommand{\spin}{\text{spin}}

\newcommand{\sgn}{\operatorname{sgn}}

\newcommand{\lk}{\operatorname{lk}}

\newcommand{\Pin}{\operatorname{Pin}} 

\newcommand{\arf}{\operatorname{Arf}} 

\newtheorem{theorem}{Theorem}
\newtheorem{proposition}{Proposition}

\theoremstyle{definition}

\newtheorem{remark}{Remark}

\newtheorem{example}{Example}

\begin{document}
\title[ ]
{10/8-inequality for a spin 4-orbifold with boundary and its application 
}

\author[]{Masaaki Ue}

\address{
Kyoto University, Kyoto, 606-8502, Japan}
\email{ue.masaaki.87m@st.kyoto-u.ac.jp}

\begin{abstract}
As an application of 10/8 inequality for a spin 4-orbifold with boundary, we give criteria of homology 3-spheres represented by certain framed links for having infinite order in the 
homology cobordism group of homology 3-spheres. Such a class of homology 3-spheres includes certain plumbed 3-manifolds, 1/odd-surgeries along certain knots with many twists, 
for example, certain Montesinos knots, $(p,q)$ cable knots of any knot with some large $p$, $q$, etc.
\end{abstract}

\maketitle

Let $Y$ be a rational homology 3-sphere with spin structure $\mathfrak s$. In \cite{U5} we describe a $10/8$-type inequality for a spin 4-orbifold bounded by $(Y, \mathfrak s )$. In this paper we apply it to give criteria of homology 3-spheres represented by certain framed links for having infinite order in the homology cobordism group $\Theta^3_H$ of homology 3-spheres. 
Such a class of homology 3-spheres includes certain plumbed homology 3-spheres and 1/odd-surgeries along certain knots with sufficiently many twists. In $\S$2 we give more explicit examples.

\section{$10/8$ type inequality for a spin 4-orbifold whose boundary is a rational homology 3-sphere and main theorem} 
Let $(X, \mathfrak s_X )$ be a compact spin 4-orbifold bounded by a rational homology 3-sphere $Y$ with $\spin$ structure $\mathfrak s$. In \cite{U5}, 
we describe a following inequality concerning the $\kappa$-invariant $\kappa (Y, \mathfrak s )$ of $(Y,  \mathfrak s )$ (\cite{Mano4}), which is deduced from the theory of Seiberg-Witten Floer stable homotopy type \cite{Mano1} and the data of $X$.  We omit to describe the unique spin structure
of $Y$ if $Y$ is an integral homology 3-sphere.

\footnote{
2020 Mathematics Subject Classification. Primary 57K31; 57K41, 57R58
}

\begin{theorem}\cite{U5}\label{estimate}
Let $(Y, \mathfrak s )$ be a rational homology 3-sphere with $\spin$ structure, and $( X, \mathfrak s_X )$ be a compact connected $\spin$ 4-orbifold  with $b_1 (X)=0$ bounded 
by $(Y, \mathfrak s )$. We assume that  the set of the singular points of $X$ consists of isolated points $\{ p_i  \ | \ i=1 ,\dots , n \}$only.
 Then a regular neighborhood of $p_i$ is a cone $cS_i$ over some spherical 3-manifold $S_i$ with $\spin$ structure $\mathfrak s_i$ induced from $\mathfrak s_X$.
Let $b_2^+ (X)$ $($resp. $b_2^- (X)$$)$be the dimension of the maximal positive $($resp.negative$)$ definite subspace of the intersection form of $X$, and $\sigma (X) = b_2^+ (X) -b_2^-(X)$ 
be the signature of $X$. 
\begin{enumerate}
\item 
\[
\kappa (Y,  \mathfrak s ) \equiv -\frac{\sigma (X)}{8} -\sum_{i=1}^n \overline{\mu } (S_i , \mathfrak s_i ) \pmod{2\bz}.
\]
\item 
\begin{align*} 
-b_2^+(X) -\frac{\sigma (X)}{8}  &-\sum_{i=1}^n \overline{\mu} (S_i , \mathfrak s_i ) -\epsilon^+ \le \kappa (Y, \mathfrak s ) \\
&\le 
b_2^- (X) -\frac{\sigma (X) }{8}  -\sum_{i=1}^n  \overline{\mu} (S_i , \mathfrak s_i )+\epsilon^-,  
\end{align*}
where
\[
\epsilon^- = \begin{cases} 0 & (b_2^-  (X)\equiv 0 \pmod 2 ) \\
 1 & (b_2^- (X)\equiv 1 \pmod 2 ) 
 \end{cases} , \ 
 \epsilon^+ = 
 \begin{cases} 
 0 & ( b_2^+ (X)=0 ) \\ -1 & ( b_2^+ (X)\equiv 1 \pmod 2 ) \\
 -2 & (b_2^+ (X)\equiv 0 \pmod 2 , b_2^+ (X)>0 ) 
 \end{cases} 
 \]
 and $\overline{\mu} (S_i , \mathfrak s_i )$ is the Neumann-Siebenmann invariant of $(S_i , \mathfrak s _i )$ \cite{N}, \cite{Si}. 
 \end{enumerate}
 
 \end{theorem}

 Theorem \ref{estimate} is a reformulation of the following $10/8$-type theorem by Manolescu. 
 
 \begin{theorem}\label{ineq}\cite{Mano4} 

Let $(Y_i , \mathfrak s_i )$ $(i=0,1 )$ be rational homology 3-spheres with $\spin$ structures and 
$(W , \mathfrak s_W )$ be a $\spin$ cobordism with $b_1 (W) =0$ from $(Y_0 , \mathfrak s_0 )$ to $(Y_1 , \mathfrak s_1 )$. 
Then 
\[
-b_2^- (W) +\frac{\sigma (W)}{8} -1 \le \kappa (Y_0, \mathfrak s_0 ) -\kappa (Y_1 , \mathfrak s_1 ) \le b_2^+ (W) +\frac{\sigma (W)}{8} +1 . 
\]
The first term of the above inequality can be replaced by $-b_2^- (W) +\frac{\sigma (W)}{8}$ if $b_2^- (W)$ is even, and 
the last term can be replaced by $b_2^+ (W) +\frac{\sigma (W)}{8}$ if $b_2^+ (W)$ is even. 

Furthermore if $(Y_0 , \mathfrak s_0 )$ is Floer $K_G$ split, 
then the second inequality is improved as follows. 
\begin{align*}
&\kappa (Y_0, \mathfrak s_0 ) -\kappa  (Y_1 , \mathfrak s_1 )   \\
&\le \begin{cases} 
b_2^+ (W) +\frac{\sigma (W)}{8} -1 &\quad \text{if $b_2^+ (W)$ is odd}, \\
b_2^+ (W) +\frac{\sigma (W)}{8} -2 & \quad \text{if $b_2^+ (W)$ is even and $b_2^+ (W) >0$. } 
\end{cases}
\end{align*}
\end{theorem}

 It is not easy to determine the value of $\kappa (Y, \mathfrak s )$ in general, but we can apply the above inequality to show that certain homology 3-spheres have infinte order in 
 $\Theta^3_H$. 
First we recall that if a closed oriented 3-manifold $Y$ is represented by a framed link $\mathcal L = \cup_i L_i$ in $S^3$ (with integer framings), then $Y$ bounds a compact 4-manifold 
$X_{\mathcal L}$, which is the union of the 4-ball $B^4$ and 2-handles $h_i$ attached along $L_i$. Furthermore a spin structure $\mathfrak s$ of $Y$ corresponds to a choice of a 
characteristic sublink $\mathcal L_0$ of $\mathcal L$, i.e., a sublink that satisfies 
\[
\cup_{L_j \in \mathcal L_0} \lk (L_j , L_i ) \equiv \lk (L_i , L_i ) \pmod 2 
\]
for each $L_i$. Here $\lk$ is the linking number and $\lk (L_i , L_i )$ means the framing of $L_i$ (this condition does not depend on the orientation of $\mathcal L$). 
Let $x_i$ be the second homology class of $X_{\mathcal L}$ represented by a 2-cycle, which is the union of the core of $h_i$ and a surface in $D^4$ bounded by $L_i$. Then the Poincare dual of the sum of $x_i$ corresponding to $\mathcal L_0$ $\pmod 2$ is the obstruction to extending $\mathfrak s$ to the spin structure of $X_{\mathcal L}$. 
We note that if $\mathcal L_0$ consists of more than one component, then we can connect them by sliding one component over the others to obtain 
a new framed link $\mathcal L ' $ with $X_{\mathcal L} \cong X_{\mathcal L '}$ so that the characteristic sublink of $\mathcal L'$ consists of a single component. 

\begin{theorem}\label{criterion}
Let $Y$ be a rational homology 3-sphere with spin structure $\mathfrak s$, which is represented by a framed link $\mathcal L$. Suppose that the characteristic sublink of 
$\mathcal L$ corresponding to $\mathfrak s$ consists of a single component $K$ with framing $n\ne 0$. 
Then there exists a positive integer $N$ such that if $|n | >N$, a connected sum $kY$ of $k$ copies of $Y$ for any $k$ does not bound a rational homology 4-ball with 
spin structure that is the extension of the spin strucuture of $kY$ induced by $\mathfrak s$. In particular if $Y$ is an integral homology 3-sphere, then the homology cobordism class  of 
$Y$ has infinite order in $\Theta^3_H$. Here $N$ is determined by $b_2 (X_{\mathcal L})$, $\sigma (X_{\mathcal L} )$, the genus $g(K)$ of $K$ and the 
Arf invariant $\arf (K)$ of $K$. 

More explicitly, suppose that $n$ satisfies one of the following conditions;
\begin{align*} 
&(c1) \quad n>4b_2 (X_{\mathcal L} ) +5\sigma (X_{\mathcal L}) +8g(K) -4(1+\sgn n )
\\
&(c2) \quad 
n = 4b_2 (X_{\mathcal L} ) +5\sigma (X_{\mathcal L}) +8g(K) -4(1+\sgn n ) > \sigma (X_{\mathcal L} )
\\
&(c3) \quad 
 n< -4b _2 (X_{\mathcal L} ) +5\sigma (X_{\mathcal L}) -8g(K) +4(1-\sgn n )
\end{align*}
 if $\arf (K) =0$, while if $\arf (K) \ne 0$, suppose that $n+8\epsilon$ satisfies one of the following conditions;
 \begin{align*}
&(c1 ' ) \quad  n+8\epsilon > 4b_2 (X_{\mathcal L} ) +5\sigma (X_{\mathcal L} ) +8g(K) +12 -4\sgn (n+8\epsilon ) 
 \\
&(c2 ') \quad 
n+8\epsilon = 4b_2 (X_{\mathcal L} ) +5\sigma (X_{\mathcal L} ) +8g(K) +12 -4\sgn (n+8\epsilon ) > \sigma (X_{\mathcal L} )
\\
&(c3 ') \quad  
n+8\epsilon< -4b_2 (X_{\mathcal L} ) +5\sigma ({\mathcal L} ) -8g(K) -12 -4\sgn (n+8\epsilon ).
\end{align*}

Here $\epsilon =\pm 1$ and we can choose any of $\epsilon$ if $n+8\epsilon \ne 0$. 
 
 Then we have the above result. 
\end{theorem}

To prove Theorem \ref{criterion}, we first follow Kaplan's method. In \cite{Kap} Kaplan constructed a framed link representing a compact spin (framed) 4-manifold bounded by a given closed oriented 3-manifold with spin structure (with almost framing). We  start with a rational homology 3-sphere with spin structure $(Y, \mathfrak s )$ represented by a framed link 
$\mathcal L$ whose characteristic sublink consists of a single knot $K$ with framing $n\ne 0$, and construct a spin 4-orbifold bounded by $(Y, \mathfrak s )$ in place of a spin 4-manifold. 
At the first stage, we construct a framed link $\mathcal L '$ representing $Y$ whose characteristic sublink corresponding to $\mathfrak s$ is a single unknot. To do so, we adopt the 
procedure in \cite{Kap} as follows. 
First we assume that $\arf (K) =0$. If $g=g(K)>0$, $K$ bounds a Seifert surface $S$ in $S^3$ and simple closed curves $a_i$, $b_i$ $(i=1,\dots , g )$ which form a 
symplectic basis of $H_1 (S, \bz )$. We can assume that $S$ is the union of a 2-disk and $2g$ bands whose cores are the parts of $a_i$, and $b_i$. It is known that if $\arf (K) =0$, we can choose $S$ and $a_i$, $b_i$ so that $\lk (a_i , a_i^+ )$ is even for all $i$. Here $a_i^+$ is the push-off of $a_i$ along a normal vector field to $S$. 

\begin{figure}
\begin{center}
\scalebox{.7}{\includegraphics[ ]{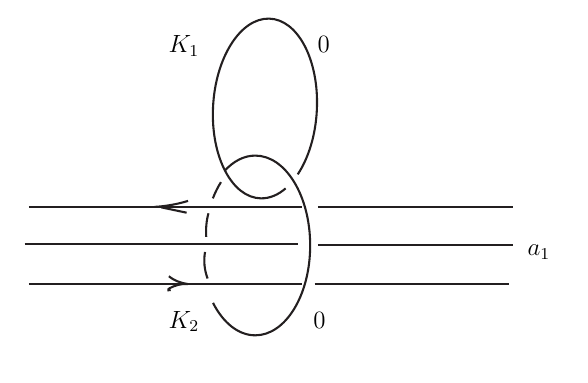}}
\caption{}
\end{center}
\end{figure}

\begin{figure}
\begin{center}
\scalebox{.6}{\includegraphics[]{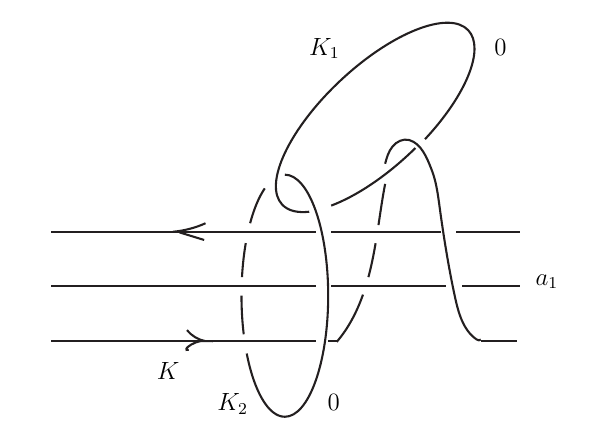}}
\caption{}
\end{center}
\end{figure}

\begin{itemize}
\item 
Add  a Hopf link $K_1 \cup K_2$ with $0$ framings to $\mathcal L$ (replace $X_{\mathcal L}$ by $X_{\mathcal L} \sharp S^2 \times S^2 $). Deform the band in $S$ by handle slides over $K_1$ so that the band links $K_2$ as in Figure 1 (Diagram 4.3 in \cite{Kap}). 
Repeat crossing changes of $a_1$ by sliding the bands containing  a part of $a_1$ over $K_2$. Then we deform $S$ so that $a_1$ is the unknot. 
\item Next we slide $K$ over $K_2$ as in Figure 2 (Diagram 4.4 in \cite{Kap}) two times so that the bands containing $a_1$ are twisted two times in the right or left direction. Since 
$\lk (K, K_2 ) =0$ and the linking number of $K_1$ and a new $K$ is $0 \pmod 2$, $K$ is still characteristic and its framing remains unchanged. 
Repeating this process, we can deform $S$ so that $\lk (a_1 , a_1^+ ) =0$. 
\item Since $a_1$ is unknotted and $\lk (a_1 , a_1^+ )$ is zero, $a_1$ bounds an embedded 2-disk $D$ such that the interior of the collar neighborhood of $\partial D$ does not intersect $S$. If $b_1$, $a_i$, $b_i$ $(i\ge 2)$i intersect the interior of $D$ (we can assume that the bands containing the parts of these curves intersect $D$), we can remove them by 
handle slides over $K_2$, and $D$ can be deformed so that its interior does not intersect $S$. Then by removing the neighborhood of $a_1 \subset S$ and attaching two copies of $D$, 
$K$ bounds a Seifert surface of genus $g-1$, which contains $a_i$, $b_i$ $(i\ge 2)$. Therefore repeating the above procedure $g$ times we obtain a new framed link 
$\mathcal L '$ with $X_{\mathcal L ' } \cong X_{\mathcal L }\sharp gS^2 \times S^2$ whose characteristic sublink consists of a single unknot with framing $n$. 
\end{itemize}

If $\arf (K) \ne 0$, add a Hopf link $K_1 \cup K_2$ with framings $0$ and $1$ (resp.  $0$ and $-1$) and slide $K_2$ over $K_1$ two times to obtain a right-handed trefoil knot with 
framing $8$ 
(resp.a left-handed trefoil knot with framing $-8$) (Diagram 4.6 in \cite{Kap}). Then sliding $K$ over the trefoil knot we obtain a new characteristic knot $K '$ with framing $n+8$ (resp.$n-8$). 
Here $g(K ' ) = g+1$, and $\arf ( K ') =0$ since the Arf invariant of the trefoil knot is $1$. Hence we can apply the above procedure to $K ' $ to obtain a 
new framed link $\mathcal L '$ with a characteristic unknot with framing $n\pm 8$, where $b_2 (X_{\mathcal L ' } ) = b_2 (X_{\mathcal L}) +2g +4$ and 
$\sigma (X_{\mathcal L ' }) = \sigma (X_{\mathcal L})$. 

In either case we have a framed link $\mathcal L ' $ with $\partial X_{\mathcal L ' }  =Y$ whose characteristic sublink corresponding to $\mathfrak s$ is the unknot 
$\widehat K$ with framing $n'$, where $n' =n$ (if $\arf (K)=0$), or $n' = n+8\epsilon$ (if $\arf (K) \ne 0 $). 
In the second case we choose $\epsilon$ so that $n+8\epsilon  \ne 0$. The union of the core of the two handle attached along 
$\widehat K$ and an embedded 2-disk bounded by $\widehat K$ in the 4-ball $D^4$ forms an embedded two sphere $S$ with $S \cdot S =n'$. 
Let $N$ be the tubular neighborhood of $S$. 
Then $\mathfrak s$ extends to the spin structure over $X_{\mathcal L ' } \setminus \mathrm{Int} N$ which induces the spin structure $\mathfrak s_0$ over the lens space 
$L( n' , -1 ) = \partial N$. Since $\mathfrak s_0$ extends to the spin structure on the cone $cL(n',-1 )$ over $L(n',-1 )$ (\cite{FFU} \cite{S}), we can replace $N$ by $cL(n', -1)$ to obtain a  
4-orbifold with spin structure $(\widehat X, \widehat{\mathfrak s})$ with $\partial (\widehat X , \widehat{\mathfrak s } ) = ( Y, \mathfrak s )$. 
Then $b_2 (\widehat X ) = b_2 (X_{\mathcal L ' } ) -1$, $\sigma (\widehat X ) = \sigma (X_{\mathcal L '} )-\sgn n' $, and 
$\overline{\mu}(L(n' , -1 ), \mathfrak s_0 ) =  (\sgn n'  -n' )/8$.
We note that $kY$ bounds a boundary connected sum of $k$ copies of $\widehat X$, which is denoted by $k\widehat X$. Hence by applying Theorem \ref{estimate} to $k\widehat X$, we have 

\[
\kappa (kY , \mathfrak s ) \equiv k (-\frac {\sigma (\widehat X )}{8} -\overline{\mu} (L( n ' , -1 ) , \mathfrak s_0)) =\frac k8 (-\sigma (X_{\mathcal L})  + n '  ) \pmod 2 
\]

where the spin structure on $kY$ induced from $\mathfrak s$ is denoted by the same symbol, and the following inequality.
\[ 
k (-b_2^+ (\widehat X) +\frac 18 (-\sigma (X_{\mathcal L}) +n ' )) - \epsilon^+ \le \kappa (kY , \mathfrak s )  
\le k (b_2^- (\widehat X ) +\frac 18 (-\sigma (X_{\mathcal L} ) +n' )) + \epsilon^- . 
\]
Here $\epsilon^+$ is the correction term with respect to $b_2^+ (k\widehat X ) =k b_2^+ (\widehat X )$ and $\epsilon^-$ is the correction term with respect to 
$b_2^- (k\widehat X ) = kb_2^- (\widehat X )$. 
If $-b_2^+ (\widehat X) +\frac 18 (-\sigma (X_{\mathcal L}) +n ' ) >0$, 
or $b_2^+ (\widehat X ) =\frac 18 (-\sigma (X_{\mathcal L})+n ' ) >0$, 
the left hand side of the above inequality is positive (in the second case $\epsilon^+ = -1$ or $-2$), while if 
$b_2^- (\widehat X ) +\frac 18 (-\sigma (X_{\mathcal L} ) +n' ) <0$, the right hand side is negative (for sufficiently large $k$). In either case $\kappa (kY , \mathfrak s ) \ne 0$. 
By rewriting these inequalities by the data of $\mathcal L$ and $K$, we obtain $(c1)\sim (c3)$, $(c1' ) \sim (c3 ')$
in Theorem \ref{criterion}. This proves Theorem \ref{criterion} due to the spin rational homology cobordism invariance of $\kappa$ invariant. 

\begin{remark}\label{family}
Let $\mathcal L$ be a framed link whose characteristic sublink is a single knot $K$ as in the assumption of Theorem \ref{criterion}. 
Suppose that there is a monotone increasing sequence of positive integers $\{ n_i \ | \ i\in \bz_{>0} \}$ and 
we have a family of rational homology 3-spheres with spin structures $(Y_i, \mathfrak s_i)$ $(i\in \bz_{>0})$ represented by framed links $\mathcal L_i$, where each $\mathcal L_i$ 
has the same link type and the same characteristic knot $K$ as $\mathcal L$ but the framing of $K$ is $n_i$. Let $\widehat X_i$ be the 
spin 4-orbifold constructed from $X_{\mathcal L_i}$ as in the proof of Theorem \ref{criterion}. Since $b_2 (X_{\mathcal L_i})$ is independent of $n_i$, 
$|\sigma (X_{\mathcal L_i})|$ and $b_2^{\pm} (\widehat X_i)$ are uniformly bounded. Thus by passing to a subsequence if necessary, we can choose $\{ n_i \} $ so that 
$-b_2^+ (\widehat X_i) +\frac 18 (-\sigma (X_{\mathcal L_i}) +n_i ' ) >0$ (where $n_i ' =n_i$ or $n_i \pm 8 $ according to the value of $\arf (K)$), and 
\[
b_2^- (\widehat X_i) +\frac 18 (-\sigma (X_{\mathcal L_i }) +n_i  ' )  +1< -b_2^+ (\widehat X_{i+1}) +\frac 18 (-\sigma (X_{\mathcal L_{i+1}}) +n_{i+1} ' )
\]
for every $i$. Thus $\kappa (Y_i , \mathfrak s_i )$ are mutually distinct and tend to $\infty$. If all $Y_i$ are integral homology 3-spheres (In \S 2 we will give such examples), we have an 
infinite family of homology 3-spheres whose homology cobordism classes are mutually distinct. 
However we do not know whether the $\kappa$ invariant of a linear combination of some $Y_i$'s is zero or not, if there are both positive and negative coefficients in the linear combination. 
\end{remark}

\begin{figure}
\begin{center}
\scalebox{.6}{\includegraphics[]{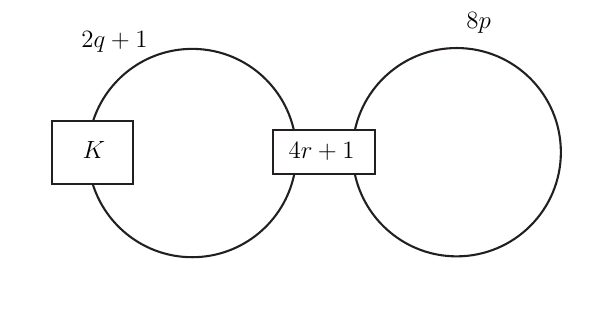}}
\caption{}
\end{center}
\end{figure}
\begin{remark}\label{arbitrary}
The conditions in Theorem \ref{criterion} only consist of the data of the characteristic knot and the linking matrix of $\mathcal L$. 
Hence if there exists an example that satisfies the above condition, we obtain infinitely many examples by replacing the components other than 
the characteristic knot by arbitrary knots. 
\end{remark}

\begin{example}\label{2component}
Let $Y$ be a 3-manifold represented by a 2 component framed link $\mathcal L$ in Figure 3, where $p,q,r \in \bz$, $K$ is an arbitrary knot, and 
the characteristic knot of $\mathcal L$ is the unknot whose framing is $8p$. If $p(2q+1) =r(2r+1)$ then the determinant of the linking matrix 
of $\mathcal L$ is $-1$ and hence $Y$ is an integral homology 3-sphere. In this case $Y$ satisfies the condition $(c1)$ if $p>0$, the condition 
$(c3)$ if $p<0$, and hence $Y$ has infinite order in $\Theta^3_H$. On the other hand if $p=0$, then $r=0$ and $Y=S^3$. In this case $Y$ 
does not satisfy any of the conditions in Theorem \ref{criterion}. 
\end{example}

\begin{remark}\label{compare}
We compare Theorem \ref{estimate} with the inequality deduced from the orbifold 10/8 theorem \cite{FF}, which states that 
\[
\ind_{\bc} \dir_Z =0 \quad \text{or} \quad 1-b_2^- (Z) \le \ind_{\bc} \dir_Z  \le b_2^+ (Z) -1
\]
for a closed 4-orbifold with spin structure $(Z, \mathfrak s_Z )$. Here $\dir_Z$ is the spin Dirac operator of $Z$.

Let $(X, \mathfrak s_X )$ be a compact spin 4-orbifold bounded by a rational homology 3-sphere with spin structure $(Y, \mathfrak s )$ as in Theorem 1. Choose a compact 4-manifold with spin structure $(W, \mathfrak s_W)$ bounded by $(Y, \mathfrak s )$ 
and let $Z=X\cup -W$ with spin structure $\mathfrak s_Z$ induced from $\mathfrak s_X$ and $\mathfrak s_W$. 
Then the Fukumoto-Furuta invariant for $(Y,X, \mathfrak s_X )$ is defined as 
\[
w(Y,X, \mathfrak s_X ) = -\ind_{\bc} \dir_Z  + \frac{\sigma (W)}{8}. 
\]
The sign convention here is opposite to the original definition in \cite{FF}. The orbifold index theorem shows that $w(Y, X, \mathfrak s_X )$ is also represented as 
\[
w(Y, X, \mathfrak s_X ) = \frac{\sigma (X)}{8} +\sum_i \overline{\mu} (S_i , \mathfrak s_i ),
\]
where $(S_i , \mathfrak s_i )$ is the spherical 3-manifold, which is the link of the singularity of $X$ with spin structure induced from $\mathfrak s_X$ as in Theorem 1 
\cite{U4} \cite{U5}. 
Then by the orbifold 10/8-theorem and the above formula, we obtain the following two inequalities. 

\[
-b_2^+ (W) -\frac{\sigma (W)}{8} \le b_2^- ( X) - \frac{\sigma(X)}{8} -\sum_i \overline{\mu} (S_i, \mathfrak s_i ) -1, 
\]
\[
-b_2^+ (X) -\frac{\sigma (X)}{8}- \sum_I \overline{\mu} (S_i , \mathfrak s_i ) \le b_2^- (W) -\frac{\sigma (W)}{8}  -1. 
\]

On the other hand Theorem \ref{estimate} shows that 
\[
-b_2^+ (W) -\frac{\sigma (W)}{8} - \epsilon_W^+ \le \kappa (Y, \mathfrak s ) \le b_2^- (X) -\frac{\sigma(X)}{8} -\sum_i \overline{\mu} (S_i , \mathfrak s_i ) +\epsilon^-, 
\]
\[
-b_2^+ (X) -\frac{\sigma (X)}{8} -\sum_i \overline{\mu} (S_i , \mathfrak s_i )-\epsilon^+  \le \kappa (Y, \mathfrak s ) \le b_2^- (W) -\frac{\sigma (W) }{8} +\epsilon_W^- ,
\]
where 
\begin{align*}
&\epsilon_W^+ =
\begin{cases} -1 &\quad (b_2^+ (W) \equiv 1 \pmod 2 ) \\
-2 &\quad (b_2^+(W) \equiv 0 \pmod 2 \ \text{and} \ b_2^+ (W) >0) \\
0 &\quad (b_2 ^+ (W) =0),
\end{cases} 
 \\
 &
 \epsilon_W^- = 
 \begin{cases} 
 1 &\quad (b_2^- (W) \equiv 1 \pmod 2) \\
 0 &\quad (b_2^- (W) \equiv 0 \pmod 2 ) ,
 \end{cases} 
 \end{align*}
 
 and $\epsilon^{\pm}$ are the correction terms that appear in Theorem \ref{estimate}. The above two inequalities with respect to $W$ are deduced from Manolescu's Theorem 
\cite{Mano4} due to the 
 fact that $S^3$ with its unique spin structure is Floer $K_G$ split. 
 Thus these inequalities are slightly different from those deduced from the orbifold 10/8 theorem. 

\end{remark}

\section{Application to plumbed homology 3-spheres and Dehn surgeries along certain knots.} 

First we consider a homology 3-sphere $Y(\Gamma )$ that bounds a plumbing $X(\Gamma )$ associated with a weighted graph $\Gamma$. Let $V =\{ v_i \}$ be the set of vertices of 
$\Gamma$ and let $m(v_i )$ be the weight of $v_i$. Here $X ( \Gamma )$ is constructed from $D^2$ bundles over $S^2$ with Euler class $m(v_i )$ for $v_i$ by patching 
them together if the corresponding vertices are connected by an edge. Then the zero section of the $D^2$ bundle over $S^2$ corresponding to $v_i$ forms a basis of 
$H_2 (X (\Gamma ), \bz )$, which we denote by $x_i$. If $Y( \Gamma )$ is a homology 3-sphere, there exists a unique characteristic element $w(\Gamma )$ of 
$H_2 (X(\Gamma ), \bz )$ corresponding to the unique spin structure of $Y(\Gamma )$ such that $w(\Gamma ) \cdot x_i \equiv  x_i \cdot x_i \pmod 2$ and $w(\Gamma ) =\sum_i \epsilon_i x_i $, where 
$\epsilon_i = 1$ or $0$. Then the Neumann-Siebenmann invariant $\overline{\mu} (Y(\Gamma ) )$ is defined as $(\sigma (X(\Gamma ))   -w(\Gamma )^2 )/8$. 

\begin{theorem}\label{plumbed}
Let $Y(\Gamma )$ be a plumbed homology 3-sphere which bounds a negative definite plumbing $X(\Gamma )$. If $\overline{\mu } (Y (\Gamma ) ) <0$ or 
$\overline{\mu} (Y (\Gamma )) >b_2 (X (\Gamma)  ) -1$, then $Y(\Gamma )$ has infinite order in $\Theta^3_H$. 
\end{theorem}
 
\begin{proof}
We note that $X(\Gamma )$ is represented by a framed link $\mathcal L  =\sum L_i$ such that $L_i$ corresponding to $v_i$ is the unknot with framing $m(v_i )$.
Suppose that $w(\Gamma ) =0$ for $Y(\Gamma )$. Then the intersection form of $X(\Gamma )$ is negative definite and of even type. If $kY(\Gamma)$ bounds a homology 4-ball for some 
$k>0$, then the union of the boundary connected sum of $k$ copies of $X(\Gamma )$ and the homology 4-ball bounded by $kY(\Gamma)$ forms a closed smooth negative definite 4-manifold of even type, which contradicts Donaldson's theorem. Hence $Y(\Gamma )$ has infinite order in $\Theta^3_H$. Suppose that $w(\Gamma ) \ne 0$. 
In this case the subset $V_0 =\{ v_i \in V \ | \ \epsilon_i =1 \}$ of $V$ contains no adjacent vertices since $\Gamma$ is a tree. Hence the 
characteristic sublink $\mathcal L_0$ of $\mathcal L$ consists of  mutually separated unknots corresponding to the vertices in $V_0$. It follows that if 
we connect the components of $\mathcal L_0$ by sliding one component over the others, we have a characteristic unknot with framing 
$w(\Gamma )^2 = \sum_{v_i  \in V_0 } x_i^2 =\sum_{v_i \in V_0 } m(v_i  ) <0$. 
Thus Theorem \ref{criterion} shows that if 
\[
0> w(\Gamma )^2 > 4b_2 (X(\Gamma )) +5\sigma (X(\Gamma )) -4(1+\sgn (w(\Gamma )^2  )) =\sigma (X(\Gamma ))
\]
 or 
 \begin{align*}
 w(\Gamma)^2 &<-4b_2 (X(\Gamma )) +5\sigma (X(\Gamma )) +4(1-\sgn (w(\Gamma )^2 )) =9\sigma (X(\Gamma )) +8 \\
 &= \sigma (X(\Gamma )) -8(b_2 (X (\Gamma ) -1 ),
 \end{align*}
 then $Y(\Gamma )$ has infinite order in $\Theta^3_H$. This proves the required result. 
 \end{proof}
 
 \begin{remark}\label{AR}
 Dai \cite{Dai} and Dai-Manolescu \cite{DM} proved the homology cobordism invariance of the $\overline{\mu}$ invariants for homology 3-spheres that bound 
 almost rational (AR) negative definite plumbings in the sense of N\'emethi \cite{Ne1}. If a plumbed homology 3-sphere $Y$ is AR, $\kappa (Y) = -\overline{\mu} (Y)$ or $-\overline{\mu} (Y) +2$ (Dai-Sasahira-Stoffregen \cite{DSS}). 
 A class of AR plumbed homology 3-spheres includes all $Y(\Gamma )$ whose associated 
 graph $\Gamma$ has at most one bad vertex, i.e., a vertex $v$ that satisfies $|m(v)|  < \delta (v)$. Here $\delta (v)$ is the number of edges whose ends contain $v$. 
 If $\Gamma$ has two or more bad vertices $v$ with $|m(v)| \le \delta (v) -2$, then $Y(\Gamma )$ is not AR (\cite{Ne1}). 
 We do not assume that $Y(\Gamma )$ is AR in Theorem \ref{plumbed}. In this case it seems difficult to prove the homology cobordism invariance of $\overline{\mu}$ by our method. 
 \end{remark}

\begin{remark}\label{optimal}
Theorem \ref{plumbed} suggests that the condition $(c1)$ in Theorem \ref{criterion} cannot be improved except for the case $(c2)$ in general. If the equality holds in $(c1)$ 
(our case does not satisfy $(c2)$), then $Y(\Gamma )$ with 
$\overline{\mu} (Y(\Gamma )) =0$ would also have infinite order in $\Theta^3_H$. But there exists $Y(\Gamma )$ with $\overline{\mu} (Y(\Gamma )) =0$ that bounds a 
contractible 4-manifold. 
\end{remark}

Next we apply the variant of Theorem \ref{criterion}, which sometimes provides better estimates than Theorem \ref{criterion},  to homology 3-spheres obtained by Dehn surgeries along knots in $S^3$ with  special properties. 

\begin{figure}
\begin{center}
\scalebox{.6}{\includegraphics[]{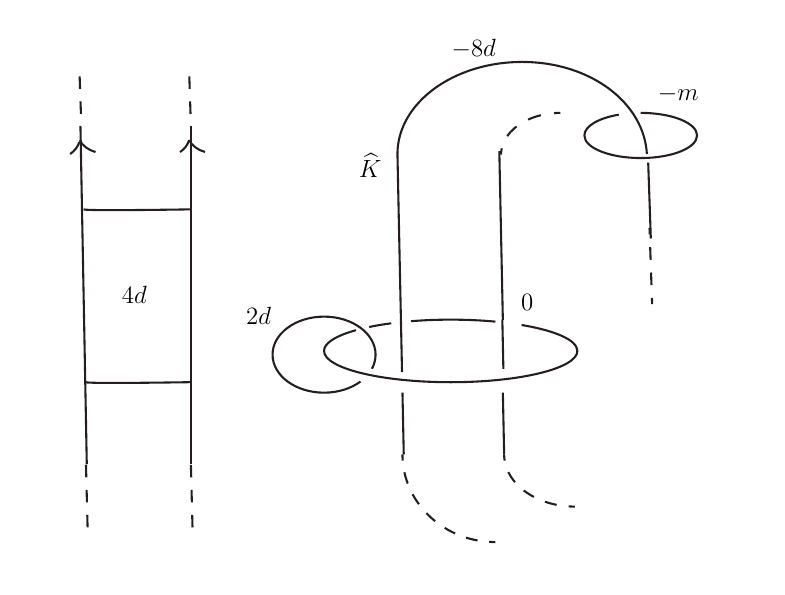}}
\caption{}
\end{center}
\end{figure}
\begin{theorem}\label{twists} 
Let $K$ be a knot in $S^3$. Suppose that $K$ contains a band with $4d$ half twists $($$2d$ full twists$)$ for $d\ne 0 $ $($right handed if $d>0$, and left handed if $d<0$$)$ as in Figure 4. We assume that if $K$ is 
oriented, the two strings of the band are oriented in the same direction. Let $\widehat K$ be a knot obtained from $K$ by removing the $4d$ half twists. Then the homology cobordism class of the homology 3-sphere $K_{1/m}$ obtained by $1/m$-surgery along $K$ for odd $m$  
has infinite order in $\Theta^3_H$ if $d$ satisfies the following conditions. 

\begin{itemize}
\item The case when $\arf  (\widehat  K ) =0$: 
\[ 
d>g(\widehat K)  \quad \text{or} \quad d\le -g(\widehat K) -1 \quad \text{if $m>0$}.
\]
\item  The case when $\arf (\widehat K ) \ne 0$: 
\[
d>g(\widehat K ) +1 \quad \text{or} \quad d \le -g(\widehat K ) -2 \quad \text{if $m>0$}.
\]


The conditions in the case when $m<0$ are obtained by replacing $d$ in the above inequalities by $-d$.. 
\end{itemize}

\end{theorem}

	\begin{proof}
	We note that $K_{1/m}$ is represented by a framed link $\mathcal L$ in Figure 4. 
	If $m$ is odd, $\widehat K$ is the characteristic knot in $\mathcal L$ corresponding to the unique spin structure on $K_{1/m}$. We also note that when $\mathcal L$ is replaced by a framed link $\mathcal L '$ which has the characteristic unknot 
$K_0$ as  in the procedure in \S 1, the meridian of $\widehat K$ with framing $-m$ remains unchanged. Hence $\mathcal L'$ contains a Hopf link consisting of $K_0$ and its meridian 
with framings $-8d$ and $-m$ if $\arf (\widehat K ) =0$, or $-8(d-\epsilon )$ $(\epsilon = \pm 1 )$ and $-m$ if $\arf (\widehat K)\ne 0$. 
It follows that $X_{\mathcal L'  }$ contains the 
plumbing $P$ represented by the above Hopf link in its interior. 
We also have an embedded 2-sphere $S'$ with $S' \cdot S' =2d$ in $\textrm{Int} X_{\mathcal L ' }$, which is the union of the core of the 2-handle attached along the unknot $c$ with 
framing $2d$ in Figure 4 and an embedded 2-disk in $D^4$ bounded by $c$, where $S'$ is disjoint from $P$. 
The boundary of $P$ is the lens space $L=L(8dm-1 , m)$ or $L=L(8(d-\epsilon )m-1 , m )$ 
(which is $S^3$ if $d=\epsilon$) according to the value of $\arf (\widehat K )$, and the boundary of 
the tubular neighborhood $N(S' )$ of $S'$ is the lens space 
$L(2d, -1 )$. 
The unique spin structure on $K_{1/m}$ extends to the spin structure on the complement of $P \cup N (S' )$, which induces the spin structures 
$\mathfrak s_0$ on $L$ and $\mathfrak s_0 ' $ on $L ( 2d, -1 )$. Hence
if we replace 
$P$ by the cone over $L$ and replace $N(S' )$ by the cone over $L(2d, -1 )$,
we obtain a spin 4-orbifold $\widehat X$. 
Let $\sigma$ (resp. $\sigma '$) be the signature of $P$ if $\arf (\widehat K )=0$ (resp. $\arf (\widehat K ) \ne 0$). 
Here they are given as follows: 
\[
 \sigma  = \begin{cases} 
 2 &\quad (\text{if $d<0$, $m<0$}) \\
 -2 &\quad  (\text{if $d>0$, $m>0$}) \\
 0 &\quad  (\text{otherwise}) 
 \end{cases} , \quad 
 \sigma ' = \begin{cases} 
   2 &\quad  (\text{if $d<\epsilon$, $m<0$}) 
   \\
   -2 &\quad   (\text{if $d>\epsilon$, $m>0$}) 
   \\
   0&\quad (\text{otherwise}). 
   \end{cases} 
   \]
Now we have $b_2 (X_{\mathcal L } ) =4$ and $\sigma (X_{\mathcal L } ) =0$ (since the intersection matrix of $X_{\mathcal L }$ is congruent to the diagonal matrix whose 
diagonal entries are $2d$, $-1/2d$, $-m$, $1/m$).  Hence
$b_2 (\widehat X ) = b_2 (X_{\mathcal L ' } ) -3 =2g(\widehat K ) +1$, $\sigma (\widehat X ) = 
\sigma (X_{\widehat L' }) -\sigma -\sgn d =-\sigma -\sgn d$, $\overline{\mu} (L, \mathfrak  s_0 ) = (\sigma +8d)/8$, and $\overline{\mu} (L(2d, -1 ) \mathfrak s_0 ' )= \sgn d /8$ if 
$\arf (\widehat K  ) =0$. 
If $\arf (\mathcal K ) \ne 0$, we have 
$b_2 (\widehat X) = b_2 (X_{\mathcal L '  } ) -3 =2g(\widehat K ) +5$, $\sigma (\widehat X ) =
\sigma (X_{\mathcal  L '} )-\sigma ' -\sgn d = -\sigma ' -\sgn d$, $\overline{\mu} (L, \mathfrak  s_0 )=(\sigma ' +8(d-\epsilon ))/8$ 
(which is $0$ if $d=\epsilon$), 
and $\overline{\mu} (L(2d, -1 ) , \mathfrak s_0 ' ) =\sgn d /8$.  
It follows that $-\sigma (\widehat X ) /8 - (\overline{\mu} (L, \mathfrak s_0 ) + \overline{\mu } (L(2d, -1 ), \mathfrak s_0 ' ) ) =-d $ if $\arf (\widehat K ) =0$, or 
$-(d-\epsilon)$ if $\arf (\widehat  K ) \ne 0 $. 
 Thus by applying the argument in the proof of Theorem 3 to $k\widehat X$, we obtain the desired result for the case when $m>0$. 
 The equalities included in the second conditions reflect $(c2)$ and $(c2')$ in Theorem \ref{criterion}.
 We note that $-K_{1/m} =K^*_{-1/m}$, where $K^*$ is the mirror image of $K$. In our setting, $K^*$ contains the band with $-4d$ half twists and the knot obtained by removing the 
 $-4d$ half twists $\widehat{K^*}$ is $(\widehat K )^*$. Since $g((\widehat K)^* ) =g(\widehat K )$ and $\arf ((\widehat K)^* ) = 
 \arf (\widehat K )$, the conditions for the case when 
 $m<0$ is deduced from the conditions for $K^*_{1/|m|}$ which is obtained by reversing the signs of $d$ and $\epsilon$ in the above inequalities for the case when $m>0$. 
  In the case when $\arf (\widehat K ) \ne 0$, we first obtain the following condition 
 \[
d-\epsilon >g(\widehat K ) +2  \quad \text{or} \quad d-\epsilon \le -g(\widehat K ) -3 \quad \text{if $m>0$}.
\]
Consequently, the case when $d=\epsilon$ is eliminated. 
By taking into account that we can choose any of $\epsilon =\pm 1$ as long as $d\ne \epsilon$, we obtain the condition as is represented above.
\end{proof}

\begin{remark}\label{uniform}
There are no uniform bounds of the number of half twists that appear in $K$ (independent of $\widehat K$) for satisfying the claim of Theorem \ref{twists}. 
For example, we note 
that for any knot $K$, the connected sum of $K$ and its mirror image with orientation reversed $-\overline K$ is a slice knot, and the $1/m$ surgery along $K\sharp 
(-\overline K)$ bounds a contractible 4-manifold. If $K$ contains the band with $4d$ half twists as in Theorem \ref{twists}, then $\widehat K$ in Theorem \ref{twists} is obtained from 
$K\sharp (-\overline K)$ by removing the $4d$ half twists, but still have the band with $-4d$ half twists. Theorem \ref{twists} suggests that the $1/m$ surgery in such a case does not 
satisfy the above inequality conditions. 
\end{remark}

 Hereafter we will give examples of more specific knots such that the $1/m$-surgeries with $m$ odd along them have infinite order in $\Theta^3_H$. 
By the same reason as is described in the proof of Theorem \ref{twists}, we give the conditions that these knots should satisfy only in the cases when $m>0$. 
We also take account of the ambiguity of the choice of $\epsilon =\pm 1$ that appears in the proof of Theorem \ref{twists}. 
 
 \begin{figure}
\begin{center}
\scalebox{.5}{\includegraphics[]{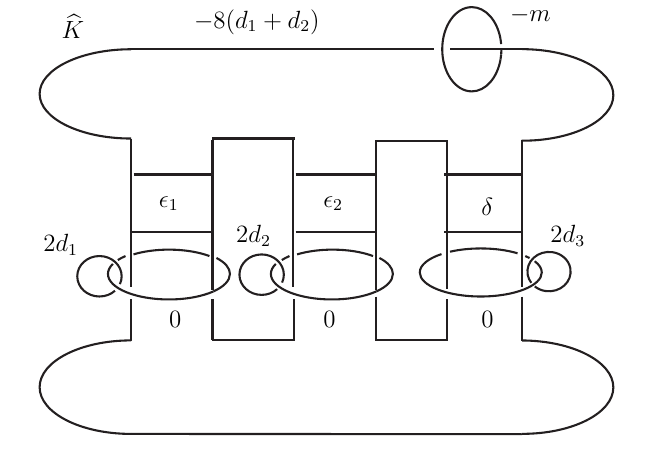}}
\caption{}
\end{center}
\end{figure}

 \begin{example}\label{pretzel}
 Let $K$ be a pretzel knot with three strands. We assume that the left two strands have 
 $4d_i +\epsilon_i$ half twists $(i=1,2 )$, and 
 the the third one has $4d_3 +\delta$ twists ($d_i \ne 0 $ for $i=1,2,3$). Here $\epsilon_i = \pm 1$ $(i=1,2 )$ and $\delta =0$ or $2$. If $K$ is oriented the two strings of the left two strands are oriented in the 
 same direction, while the two strands of the third one is oriented in the opposite direction. Then $K_{1/m}$ is represented by a framed link $\mathcal L$ in Figure 5, 
where $\widehat K$ is a knot obtained from $K$ by removing the $4d_i$ half twists $(i=1,2,3 )$ of the three strands. Then $\widehat K$ is the characteristic knot if $m$ is odd. We can see that $\widehat K$ is the unknot if 
$(\epsilon_1 , \epsilon_2 , \delta )=(\pm 1 , \pm 1 , 0 )$, $(1, -1 , 2 )$, or $(-1, 1, 2)$ and $g(\widehat K ) =1$ and $\arf (\widehat K ) =1$ otherwise 
	by considering the Seifert form of $\widehat K$. Note that $b_2 (X_{\mathcal L})=8$ and $\sigma (X_{\mathcal L } )=0$ since the intersection matrix of $X_{\mathcal L}$ is congruent to the diagonal matrix whose diagonal entries consist of $2d_i$, $-1/2d_i$ $(i=1,2,3)$, $-m$, $1/m$. 

	First suppose that $\widehat K$ is the unknot. Then $X_{\mathcal L}$ contains the plumbing $P$ corresponding the union of $\widehat K$ and its meridian with framing 
	$-8(d_1 +d_2)$ and $-m$. 	Furthermore we have embedded 2-spheres $S_i$ with $S_i \cdot S_i =2d_i$, each of which is the union of the core of the 2-handle attached along the unknot $c_i$ with framing $2d_i$ and 
	embedded 2-disks bounded by $c_i$ in $D^4$ $(i=1,2,3)$ in $X_{\mathcal L }$. Here $P'$ and the tubular neighborhoods $N(S_i )$ of $S_i$ are mutually disjoint. If $d_1 +d_2=0$, then the framed link in Figure 5 representes $S^3$. This case will be eliminated also by the argument below. The boundary of $P$ is the lens space $L =L(8m(d_1 +d_2 )-1 , m )$ ($S^3$ if $d_1 +d_2 =0$) and the boundary of $N(S_i )$ is the lens space $L(2d_i , -1 )$. 
They have spin structures $\mathfrak s_0$ and $\mathfrak s_i$ respectively,which are induced from the spin structure on $X_{\mathcal L} \setminus (P \cup \cup_i N(S_i ) )$. 
Hence replacing $P$ and $N(S_i )$  by the cones over $L$ and $L(2d_i , -1 )$, we obtain the spin 4-orbifold  $\widehat X$ bounded by $K_{1/m}$. 
Then we have 
	$b_2 (\widehat X) = b_2 (X_{\mathcal L} ) -5 =3$, $\sigma (\widehat X) = -\sigma - \sum_i \sgn d_i$, $\overline{\mu} (L, \mathfrak s_0 ) = d_1 +d_2 +\sigma /8$ (which is $0$ if $d_1 +d_2 =0$), and 
	$\overline{\mu} (L(2d_i , -1 ) , \mathfrak s_i ) = \sgn d_i /8$. Here 
	$\sigma$ is the signature of $P$, which is $-2$ if $m>0$ and $d_1 +d_2 >0$, $2$ if $m<0$ and $d_1 +d_2  <0$, and $0$ otherwise. 
	
	If $g (\widehat K )=1$ and $\arf (\widehat K ) =1$, we can replace $\mathcal L$ by a framed link $\mathcal L'$ which contains a Hopf link consisting of the characteristic unknot and its 
meridian with framings $-8(d_1 +d_2 -\epsilon )$ $(\epsilon = \pm 1 )$ and $-m$. 
Then $X_{\mathcal L '}$ contains a plumbing $P' $ represented by the above Hopf link and 
the tubular neighborhood $N(S_i )$ of the embedded 2-sphere $S_i$ with $S_i \cdot S_i =2d_i$ as in the first case. 
The lens spaces $L'=L(8m(d_1 +d_2 -\epsilon ) -1 , m ) =\partial P'$ (which is $S^3$ if $d_1 +d_2 -\epsilon =0$) 
and $L(2d_i , -1 ) =\partial N(S_i )$
admit a spin structure $\mathfrak s_0$ and $\mathfrak s_i$ respectively, which are induced from the spin structure on the complement of 
$P'$. Hence 
replacing $P'$ and $N(S_i)$ in $X_{\mathcal L ' }$ by the cones over $L'$ and $L(2d_i , -1 )$ we obtain a spin 4-orbifold $\widehat X$. Here we have 
$b_2 (\widehat X ) = b_2 ( X_{\mathcal L ' } ) -5 =9$, $\sigma (\widehat X ) = -\sigma ' -\sum_i \sgn d_i$, and 
$\overline{\mu} (L', \mathfrak s_0 ) = d_1 +d_2 -\epsilon +\sigma ' /8$ (which is $0$ if $d_1 +d_2 -\epsilon =0$)
and $\overline{\mu} (L(2d_i , -1 ) , \mathfrak s_i ) = \sgn d_i /8$, where $\sigma ' $ is the signature of $P'$. 

Then as in the proof of Theorem \ref{criterion}, 
we can see that $K_{1/m}$ for $m$ odd has infinite order in $\Theta^3_H$ under one of the following conditions: 

\begin{align*}
-&b^+_2 (\widehat X) -\frac{\sigma (\widehat X) }{8} -\overline{\mu} (L^{(')}, \mathfrak s_0 ) -\sum \overline{\mu} (L(2d_i , 
\mathfrak s_i ) >0 , \\
&b^+_2 (\widehat X) = -\frac{\sigma (\widehat X) }{8} -\overline{\mu} (L^{(')}, \mathfrak s_0 ) -\sum \overline{\mu} (L(2d_i , 
\mathfrak s_i ) >0 , \\
&b_2^- (\widehat X) -\frac{\sigma (\widehat X) }{8} -\overline{\mu} (L^{(')}, \mathfrak s_0 ) -\sum \overline{\mu} (L(2d_i  , 
\mathfrak s_i ) <0 ,
\end{align*}

from which we deduce the following conditions. (The case when $d_1 +d_2 =0$ or $d_1 +d_2 -\epsilon =0$ is eliminated in the list below.)
\begin{itemize}
\item The case when $\widehat K$ is the unknot: 
\[ 
d_1 +d_2>(1+\sum_i \sgn d_i )/2 \quad  \text{or} \quad d_1 +d_2 \le (-3 +\sum_i \sgn d_i )/2  \quad  \text{if $m>0$}.
\]
We note that $d_1 +d_2 >0$ under the first condition and $d_1 +d_2 <0$ under the second condition.

\item The case when $g(\widehat K ) =1$ and $\arf (\widehat K ) =1$: 
\[
d_1 +d_2  >(5+\sum_i \sgn d_i )/2  \quad \text{or} \quad d_1 +d_2  \le (-7 +\sum_i \sgn d_i )/2 \quad \text{if $m>0$}. 
\]
In both cases, the conditions for the cases when $m<0$ are deduced from the cases when $m>0$ by replacing $d_i$ by $-d_i$. 
\end{itemize}

If $d_i \ne 0$ $(i=1,2)$ and $d_3 =0$, the framed link representing $K_{1/m}$ is obtained from the framed link in Figure 5 by deleting the Hopf link with framings $0$ and $2d_3$. 
In this case computation similar to the above case shows that $K_{1/m}$ ($m$ odd) has infinite order in $\Theta^3_H$ under the following conditions. 

\begin{itemize}
\item The case when $\widehat K$ is the unknot:
\[
d_1 +d_2 >\sum_i \sgn d_i  /2 \quad \text{or} \quad d_1 +d_2 \le -1 +\sum_i \sgn d_i /2  \quad \text{if $m>0$}. 
\]

\item 
The case when $g(\widehat K ) =1$ and $\arf (\widehat K ) =1$:
\[
d_1 +d_2 >2+\sum_i \sgn d_i /2 \quad \text{or} \quad d_1 +d_2 \le -3 +\sum_i \sgn d_i  /2  \quad \text{if $m>0$}.
\]
\end{itemize}
The conditions for $m<0$ are also obtained from the above conditions by replacing $d_i$ by $-d_i$. 

\end{example}

\begin{remark}\label{zero} 
If $d_3 =\delta =0$, then $\widehat K$ is the unknot and $K$ is the connected sum of the torus knots $T(2, \epsilon_1 +4d_1 ) \sharp -T(2, \epsilon_2 +4d_2 )$. 
If $d_1$ and $d_2$ have opposite sign, the above condition $d_1 +d_2 >\sum_i \sgn d_i /2 =0$ cannot be improved. If $d_1 +d_2 =0$ and we further assume that 
$\epsilon_1 = -\epsilon_2$, then 
$K= T(2, \epsilon_1 +4d_1 ) \sharp -(T(2, \epsilon_1 +4d_1 )^*$ is a slice knot and hence 
$K_{1/m}$ bounds a contractible 4-manifold.

\end{remark}

In case of a 2-bridge knot $K=C(a_1 , \dots , a_n )$, the variant of Theorem \ref{criterion} can be applied under some parity conditions on $a_i$. Suppose that $a_1$ is odd and 
$a_i$ $(i\ge 2)$ is even. Then if $K$ is oriented, the two upper strings, where $a_{2i+1}$ half twists lie between them, are oriented in the same direction if $n$ is odd, and hence 
the variant of Theorem 3 is applied to the $a_{2i+1}$ twists. If $n$ is even, 
the lower two strings are oriented in the same direction. We only give the result in the first case. 

 \begin{figure}
\begin{center}
\scalebox{.6}{\includegraphics[]{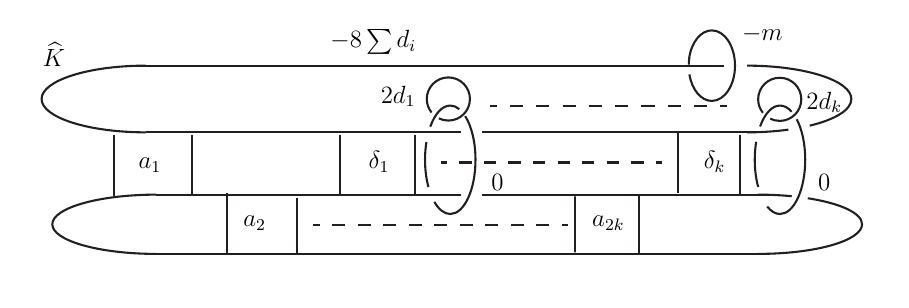}}
\caption{}
\end{center}
\end{figure}

\begin{example}\label{2bridge}
Let $K = C(a_1 , a_2 , \dots , a_{2k+1} )$ be a 2-bridge knot such that $a_1$ is odd, $a_{2i+1} = 4d_i +\delta_i$ $d_i \ne 0$ $(1\le i \le k )$, and $a_{2i}$'s are arbitrary even integers. Here we assume that $\delta_i =0$ or $2$. Then $K_{1/m}$ is represented by a framed link $\mathcal L$ in Figure 6, where $\widehat K =C(a_1 , a_2, \delta_1,   \dots , 
a_{2k}, \delta_k)$. If $m$ is odd, $\widehat K$ is a characteristic knot with framing $-8\sum d_i$. Then as in the previous case, $\mathcal L$ can be replaced by a framed link 
$\mathcal L'$ which contains a characteristic unknot $K_0$ with framing $-8\sum d_i$ if $\arf  (\widehat K ) =0$ (resp.$-8(\sum d_i -\epsilon)$ $(\epsilon =\pm 1)$ if 
$\arf (\widehat K ) \ne 0 $). 
In $X_{\mathcal L ' }$ we have a plumbing $P$ represented by the union of $K_0$ and its meridian with framing $-m$, and also the tubular neighborhoods $N(S_i)$ of embedded 2-spheres $S_i$, 
each of which is the union of the core of the 2-handle attached along the unknot $c_i$ with framing $2d_i$ in Figure 6 and an embedded 2-disk bounded by $c_i$ in $D^4$. 
Hence replacing $P$ and $N(S_i )$ by the cone over $L = \partial P$ 
(which is $S^3$ if $\sum d_i =0$ or $\sum d_i -\epsilon =0$) and the cones over $L(2d_i , -1 ) = \partial N(S_i )$ respectively, 
we obtain a spin 4-orbifold $\widehat X$ as before. We have 
 $b_2 (\widehat X ) = b_2 (X_{\mathcal L ' } ) -2-k =
 k+2g(\widehat K )$ if $\arf (\widehat K ) =0$, or $k+2(g(\widehat K) +2)$ if $\arf (\widehat K ) \ne 0$, and 
$\sigma (\widehat X ) = -\sigma ' - \sum_i \sgn d_i$. Here $\sigma '$ is the signature $P$, which is $\pm 2$ or $0$ according to the signs of 
$m$ and $\sum d_i$ or $\sum d_i -\epsilon$ as before. 
Furthermore $\overline{\mu} (L, \mathfrak s_0 ) = \sum d_i +\sigma ' /8$ if $\arf (\widehat K ) =0$ (resp. 
$\sum d_i -\epsilon +\sigma ' /8$ if $\arf (\widehat K ) \ne 0$), and $\overline{\mu} (L(2d_i , -1 ) , \mathfrak s_i ) = \sgn d_i /8$, where
$\mathfrak s_0$ and $\mathfrak s_i$ are the spin structures on $L$ and $L(2d_i , -1 )$ induced from the spin structure on the complement of $P \cup \cup _i N(S_i )$. 
 Hence we can see that 
$K_{1/m}$ with $m$ odd has infinite order in $\Theta^3_H$ under the following conditions. 

\begin{itemize}
\item The case when $\arf (\widehat K ) =0$:
\[
\sum d_i  >(k+\sum_i \sgn d_i )/2 +g(\widehat K ) -1 \ \text{or} \ \sum d_i \le (-k+\sum_i \sgn d_i)/2  -g(\widehat K )  \ \text{if $m>0$}. 
\]

\item The case when $\arf (\widehat K )\ne 0$:
\begin{align*}
&\sum d_i  > (k+\sum_i \sgn d_i )/2 +g(\widehat K )  \ \text{or} \\
&\sum d_i  \le (-k+\sum_i \sgn d_i )/2 -g(\widehat K ) -1  \ \text{if $m>0$}. 
\end{align*}
The case when $\sum d_i =0$ or $\sum d_i -\epsilon =0$ $(\epsilon =\pm 1 )$ is also eliminated in the above list. 
 The conditions for the cases when $m<0$ are deduced from the cases when $m>0$ by replacing $d_i$ by $-d_i$. 
\end{itemize} 
We note that there is a formula for computing the genera of 2-bridge knots (\cite{BZ} 12.24 Proposition).
\end{example}

The variant of Theorem \ref{criterion} also can be applied to Montesinos knots with several rational tangles if the Conway representations of these tangles satisfy certain parity conditions. We omit the 
details. 

Next we consider $(p,q)$ cable knots. 

 \begin{figure}
\begin{center}
\scalebox{.6}{\includegraphics[]{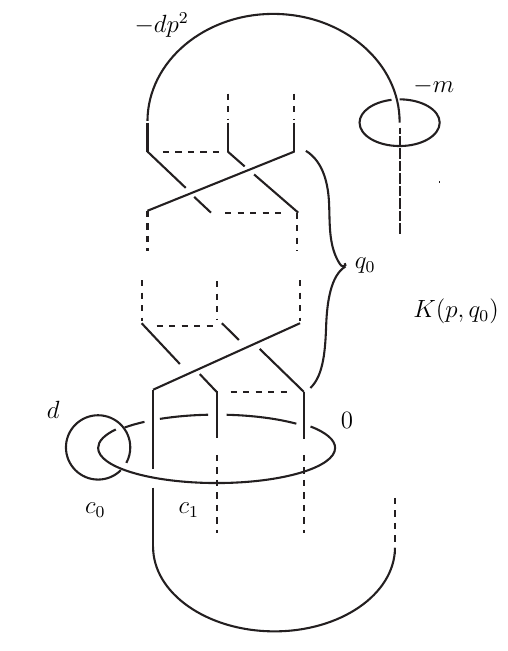}}
\caption{}
\end{center}
\end{figure}
\begin{example}\label{cable}
Let $K$ be any knot in $S^3$ and $K(p,q)$ be a $(p,q)$ cable knot of $K$. We assume that $q=q_0 +dp$ such that $d\ne 0$, $|q_0 | <p$, and $\gcd (p,q_0 ) =1$. 
Then $K(p,q)_{1/m}$ is represented by a framed link $\mathcal L$ in Figure 7.  Suppose that $m$ is odd. Then if both $p$ and $d$ are even, 
$K(p,q_0)$ is a characteristic knot with framing $-dp^2$. If $p$ is even and $d$ is odd, then the union of $K(p,q_0 )$ and the component $c_1$ with framing $0$ is a characteristic sublink of $\mathcal L$, and hence sliding $K(p,q_0 )$ over $c_1$ we obtain a new characteristic knot $\widehat K_1$ with framing $-dp^2 +2p$ and also 
$\widehat K_{-1}$ with framing $-dp^2 -2p$. If $p$ is odd, then 
the union of $K(p,q_0)$ and the component $c_0$ with framing $d$ is a characteristic sublink, and hence sliding $K(p,q_0)$ over $c_0$ we obtain a characteristic knot, which has the same knot type as $K(p,q_0 )$ but its framing becomes $-dp^2 +d$. In either case we denote the new characteristic knot by the same symbol $\widehat K$. 
Let $\delta$ be $dp^2$ if $p\equiv d \equiv 0 \pmod 2$, $dp^2 -2\epsilon p$ if $p\equiv 0$, $d\equiv 1 \pmod 2$ and $\widehat K=\widehat K_{\epsilon}$ $(\epsilon =\pm 1 )$, and $dp^2 -d $ if $p\equiv 1 \pmod 2$. 

We replace $\mathcal L$ by a framed link $\mathcal L'$ which contains a Hopf link consisting of the characteristic unknot with framing $-\delta$ if $\arf (\widehat K) =0$ 
(resp. $-\delta +8\epsilon$ $(\epsilon =\pm 1 )$ if $\arf (\widehat K ) \ne 0 $) and its meridian with framing $-m$. 
Then further replacing the pluming $P$ in $X_{\mathcal L ' }$ corresponding to the above Hopf link by the cone over the lens space $L =\partial P$, we obtain a spin 4-orbifold $\widehat X$ as before. Here $b_2 (\widehat X ) = 2(g(\widehat K )+1)$ if $\arf (\widehat K ) =0$, $2(g(\widehat K ) +3)$ if $\arf (\widehat K ) \ne 0$, $\sigma (\widehat X ) = -\sigma '$ where 
$\sigma ' $ is the signature of $P$, and $\overline{\mu} (L, \mathfrak s_0 ) = (\sigma ' +\delta )/8$ if $\arf (\widehat K ) =0$, $(\sigma ' +\delta -8\epsilon )/8$ if 
$\arf(\widehat K ) \ne 0$ where $\mathfrak s_0$ is a spin structure induced from the spin structure on the complement of $P$. 
(If $-\delta +8\epsilon =0$, then $L=S^3$ and $\overline{\mu} (L, \mathfrak s_0 ) =0$). 
We also note that if $p\equiv d\equiv 0\pmod 2$, the embedded 2-sphere $S$ in $X_{\mathcal L}$ consisting of the core of the 2-handle attached along the unknot $c_0$ with 
framing $d$ in Figure 7 and an embedded 2-disk bounded by $c_0$ in $D^4$ still exists in $X_{\mathfrak L '}$ is disjoint from $P$ since the characteristic unknot obtained by 
the procedure in \S 1 is separated from $c_0$. Hence we can further replace the tubular neighborhood $N(S)$ by the cone over $L(d, -1 ) = 
\partial N(S)$ to obtain a spin 4-orbifold 
$\widehat X$. But in this case we have the same estimates as the other cases consequently.
Applying the variant of Theorem \ref{criterion} to $\widehat X$ again, we can see that $K(p,q)_{1/m}$ with $m$ odd has infinite order in $\Theta^3_H$ under the following conditions. We can see that the case when $\delta =8\epsilon$ $(\epsilon =\pm 1)$ does not appear in the second list below. 

\begin{itemize}
\item The case when $\arf (\widehat K )=0$, 
\[
 \delta >8g(\widehat K )  \quad \text{or} \quad \delta \le -8(g(\widehat K )  +1 ) \quad \text{if $m>0$}.
 \]

 \item The case when $\arf (\widehat K) \ne 0$, 
  \[
 \delta >8(g(\widehat K ) +1) \quad \text{or} \quad \delta  \le -8(g(\widehat K ) +2) \quad \text{if $m>0$}.
\]

\end{itemize} 
Here $\delta$ is defined as 
\[
\delta =\begin{cases} 
dp^2 &\quad (p\equiv d\equiv 0 \pmod 2 \  \text{and $\widehat K=K(p,q_0 )$}), \\
dp^2 -2\epsilon p &\quad (p\equiv 0, d \equiv 1 \pmod 2 \ \text{and $\widehat K = \widehat K_{\epsilon}$)} \ (\epsilon =\pm 1 ),\\
dp^2 -d &\quad (p\equiv 1 \pmod 2, \ \text{and $\widehat K = K(p,q_0 )$}. 
\end{cases}
\]
The conditions for the cases when $m<0$ are deduced from those for the cases when $m>0$ by replacing $d$ by $-d$, $\delta$ by $-\delta$.

 \end{example}

\begin{remark}\label{restriction}
\begin{enumerate}
\item The $1/2n$-surgery along a knot $K$ in $S^3$ is represented by a framed link $\mathcal L$ consisting of $K$ with framing $0$ and its meridian with framing $-2n$. But the 
characteristic sublink of $\mathcal L$ is empty, and hence we cannot apply our method to this case. In fact, we can see that 
for any knot $K$, we have $\kappa (K_{1/2n}) =0$ if $n>0$, and $0$ or $2$ if $n<0$ (\cite{U5} Theorem 5). 

\item The statement of Theorem \ref{twists} in the case when $\arf (\widehat K)=0$ is false if $m$ is even. If $\widehat K$ is the unknot, the second condition for $m>0$ becomes 
$d\le -1$. If we add $(-4)$-twists to the unknot $\widehat K=T(2,1 )$, we obtain $K=T(2, -3)$. For $m>0$ we have $T(2,-3)_{1/m} = 
-\Sigma (2,3,6m+1 )$, which bounds a contractible 4-manifold if $m=2$ or $4$. 

\item We cannot detect the difference between $K_{1/m}$ and $K_{1/m'}$ for $m$, $m '$ odd and $m\ne m'$ because 
if the inequality in the proof of Theorem 3 is applied to $K_{1/m}$, it does not reflect the value of $m$. 
In fact, there exists a knot $K$ (for example, the knot in \cite{U5} Proposition 5), such that $\kappa (K_{1/m})$ has the same value for every $m$ which is odd and positive.
We also note that $K_{1/m}$ bounds a contractible 4-manifold for any $m\ne 0$ if $K$ is a slice knot (Proposition \ref{contractible}). But we do not know whether $K_{1/m}$ and 
$K_{1/m'}$ could be homology cobordant if $m\ne m'$ and $\kappa (K_{1/m} )=\kappa (K_{1/m '}) \ne  0$. 
\end{enumerate}
\end{remark}

\begin{remark}\label{family2}
There are several families of homology 3-spheres consisting of $1/n$-surgeries $(n\in \bz_{>0})$ along a fixed knot that generate $\bz^{\infty}$ subgroups of 
$\Theta^3_H$. Such a class of knots includes the torus knots (\cite{FS},\cite{F}), and knots that satisfy the certain conditions on the Floer theoretical invariants concerning the knots (\cite{BS1},\cite{BS2},\cite{D},\cite{DIST},\cite{NST}). We refer \cite{Sav} for more references. But we do not know what happens when $n$ is negative in general. It is also uncertain whether the examples in this section could satisfy the conditions 
on the invariants that appear in the above results. 
\end{remark}

Next we consider Dehn surgeries along slice knots. 
We recall the following fact. 

\begin{proposition}\label{contractible}
Let $K_{1/n}$ $(n\ne 0 )$ be a $1/n$-surgery along a slice knot $K$. Then $K_{1/n}$ bounds a contractible 4-manifold. 
\end{proposition}

We note that $K_{1/n}$ is represented by a framed link consisting of $K$ with framing $0$ and its meridian $m$ with framing $-n$. In the corresponding 4-manifold $X$ we have 
an embedded 2-sphere $S$ with $S\cdot S=0$, which is the union of the core of the two handle attached along $K$ and an embedded 2-disk bounded by $K$ in $D^4$. We also have an 
embedded 2-sphere $S'$ with $S' \cdot S' = -n$, which is the union of the core of the two handle attached along $m$ and an embedded 2-disk bounded by $m$ in $D^4$. 
Then the regular neighborhood of $S\cup S'$ in $X$ is a plumbing whose boundary is $S^3$, and hence replacing this plumbing with the 4-ball, we obtain a contractible 4-manifold 
bounded by $K_{1/n}$. 

On the other hand for a rational homology 3-sphere $K_m$ obtained by $m$-surgery along $K$ $(m\ne 0 )$ we have the following result. 
Consider $K_m$ as the union of the complement $S^3 \setminus \mathrm{Int} N(K)$ of the 
tubular neighborhood $N(K)$ of $K$ and a newly attached solid torus. If $m$ is even, there exists two spin structures on $K_m$: One of them extends to the unique spin structure on 
$S^3$ when it is restricted to $S^3 \setminus 
\mathrm{Int} N(K)$, which is denoted by $\mathfrak s_0$, while the other does not, which is denoted by $\mathfrak s_1$. 
If $m$ is odd, only $\mathfrak s_1$ exists. 

\begin{proposition}\label{rational}
Let $K_m$ be a rational homology 3-sphere obtained by $m$-surgery along a slice knot $K$ $(m\ne  0 )$, and $\mathfrak s_i$ $(i=0,1 )$ be the spin structures on $K_m$ as above. 
Then we have 
\[
 \kappa (K_m , \mathfrak s ) = 
 \begin{cases}
 -\sgn m /8 &\quad (m \equiv 0 \pmod 2, \ \mathfrak s = \mathfrak s_0 ), \\
- (\sgn m -m )/8&\quad (m\equiv 0 \pmod 2,  \ \mathfrak s = \mathfrak s_1 ) ,\\
 -(\sgn m -m )/8 &\quad (m\equiv 1 \pmod 2,  \ \mathfrak s=\mathfrak s_1 ). 
 \end{cases} 
 \]
 \end{proposition}
 
  \begin{proof} 
  Since $K_m$ is represented by a framed link consisting of $K$ with framing $m$, the corresponding 4-manifold $X$ contains an embedded 2-sphere $S$ with 
  $S \cdot S=m$, which is the union of the core of the 2-handle attached along $K$ and an embedded 2-disk bounded by $K$ in $D^4$. Hence replacing 
  the tubular neighborhood $N(S)$ of $S$ by the cone over the lens space 
  $L(n, -1 ) = \partial N(S)$, we obtain a spin 4-manifold $\widehat X$. Here $\mathfrak s$ extends to a spin structure on the complement of 
 the tubular neighborhood of $S$, and hence induces the spin structure $\mathfrak s'$ on $L(n,-1)$.  
If $\mathfrak s=\mathfrak s_1$, then $S$ represents a characteristic element of $X$ and hence $\overline{\mu} (L(n,-1), \mathfrak s' ) =(\sgn m -m )/8$, while if $\mathfrak s = 
\mathfrak s _0$, $\overline{\mu} ( L (n,-1 ) , \mathfrak s' ) = \sgn n /8$. We have also $b_2 (\widehat X) = \sigma (\widehat X ) =0$, and hence $b^{\pm }_2 (\widehat X  ) =0$ and 
$\epsilon^{\pm } =0$ when we apply the inequality in Theorem \ref{estimate} to $\widehat X$. Thus we obtain the desired result.  We note that $\kappa (K_{\pm 1 } , \mathfrak s_1 ) =0$, which is 
consistent with the fact that $K_{\pm 1 }$ bounds a contractible 4-manifold. But $\kappa (K_m  , \mathfrak s ) \ne 0$ If $m\ne \pm 1$, which means that 
$K_m$ does not bound a spin  rational homology 4-ball. 
\end{proof}

\begin{remark}\label{nonslice}
If $\kappa (K_{1/m} ) \ne 0$ ($m$ odd) for a knot $K$ in $S^3$, $K$ is not a slice knot. But it is uncertain to what extent this result is effective for detecting the non-sliceness of knots. 
\end{remark}

Finally we give comments on cyclic branched covers of $S^3$ branched along knots. 

If a knot $K$ in $S^3$ bounds a Seifert surface $F$, pushing $F$ into $D^4$ we obtain a properly embedded surface $\widetilde{F}$ in $D^4$. Let $X$ be a $p$-fold cyclic 
branched cover of $D^4$ branched along $\widetilde F$, which is bounded by the $p$-fold cyclic branched cover of $S^3$ branched along $K$. Akbulut-Kirby's method \cite{AK} provides 
a framed link representing $X$. But all components of the obtained framed link have even framings, and hence our method is not efficient for this case. 
On the other hand if $K$ bounds a non-orientable suface $F$ in $S^3$, we can construct a double branched cover $X$ of $D^4$ branched along $\widetilde F$, which is obtained by pushing $F$ into $D^4$ (in this case only a double branched cover can be constructed). Sometimes we can apply our method effectively to such a case. 

 \begin{figure}
\begin{center}
\scalebox{.6}{\includegraphics[]{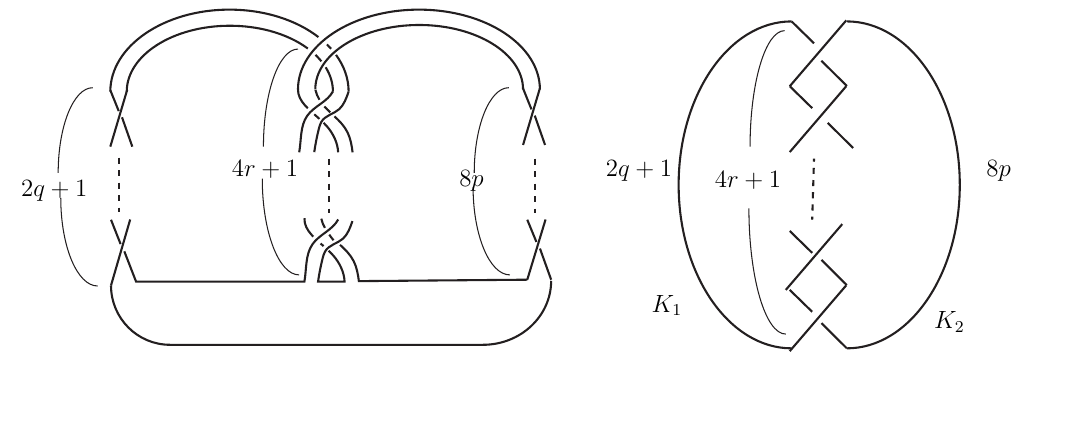}}
\caption{}
\end{center}
\end{figure}
\begin{example}\label{branched}
Let $K$ be a knot in $S^3$ bounding a non-orientable surface $F$ in Figure 8. 
Here $F$ consists of a union of a 2-disk and two bands with $2q+1$ half twists and 
$8p$ half twists respectively and they are crossed  $4r+1$ times. Let $\widetilde F$ be a surface in $D^4$ obtained by pushing $F$ into $D^4$, and let $X$ be a double branched cover of $D^4$ branched along $\widetilde F$. Then $X$ is represented by a framed link in Figure 8, which consists of two unknots $K_1$ and $K_2$ with framings $2q+1$ and $8p$ respectively such that $\lk (K_1 , K_2) = 4r+1$. 
This framed link is the same as $\mathcal L$ in the case when $K$ is the unknot that appeared in  Figure 3. Hence 
if $p(2q+1) =r(2r+1)$, the double branched cover $Y=\partial X$ of $S^3$ 
branched along $K$ is an integral homology 3-sphere, and has infinite order in $\Theta^3_H$ if $p\ne 0$ as in Example \ref{2component}. On the other hand if $p=0$, then $K$ is the unknot and $Y=S^3$. 

  \end{example}
  
  \subsection*{Tool Use} No AI tools were used in the preparation of this work.

\end{document}